\def\timestamp{%
Time-stamp: <hi-arxiv.tex: Thursday 26-10-2000 at 12:03:23 (cest)>}
\def\stripname Time-stamp: <#1 #2>{#2}
\edef\filedate{\expandafter\stripname\timestamp}

\documentclass[a4paper]{amsart}[2000/06/02]
\usepackage{eucal,amsfonts}
\newcommand\calA{\mathcal{A}}
\newcommand\calB{\mathcal{B}}
\newcommand\calC{\mathcal{C}}
\newcommand\calD{\mathcal{D}}
\newcommand\calF{\mathcal{F}}

\newcommand\calN{\mathcal{N}}
\newcommand\calT{\mathcal{T}}
\renewcommand\newsymbol[5]{%
\DeclareMathSymbol#1{#3}{\ifcase #2\or AMSa\or AMSb\fi}{"#4#5}}
\newcommand\newbbbletter[2]{%
\DeclareMathSymbol#1{0}{AMSb}{`#2}}
\let\emptyset \undefined
\let\ge       \undefined
\let\le       \undefined
\newsymbol\restr       1216
\newsymbol\mapdiagbin  124D    
\newsymbol\le          1336    
\newsymbol\ge          133E    
\newsymbol\emptyset    203F
\newsymbol\notle       230A
\newsymbol\subsetneq   2328
\newsymbol\nsubseteq   232A
\newbbbletter\Cantor   C
\newbbbletter\I        I
\let\L\undefined
\newbbbletter\L        L
\newbbbletter\N        N
\newbbbletter\psarc    P
\newbbbletter\R        R
\newcommand\cl[1]{\operatorname{cl}#1}
\newcommand\diam{\operatorname{diam}}

\newcommand\conn{\operatorname{conn}}
\newcommand\Hcube{\I^\infty}
\newcommand\preim{^{-1}}
\DeclareMathSymbol\mapdiag{1}{symbols}{"34}
\newcommand\Hyp[1]{2^{#1}}                 
\newcommand\HypX{\Hyp{X}}
\newcommand\Cs[1]{\calC(#1)}               
\newcommand\CsX{\Cs{X}}
\newcommand\Viet[1]{\langle #1\rangle}
\newcommand\wlev[1]{\mu\preim(#1)}
\let\epsilon\varepsilon
\let\meet\wedge
\let\join\vee
\let\bigmeet\bigwedge

\newcommand\0{\mathbf{0}}
\newcommand\1{\mathbf{1}}
\let\implies\rightarrow
\DeclareMathSymbol{\lor}{\mathrel}{symbols}{"5F}
\DeclareMathSymbol{\land}{\mathrel}{symbols}{"5E}
\def\nfrac#1/#2{\leavevmode\kern.1em
  \raise.5ex\hbox{\the\scriptfont0 #1}\kern-.1em
  /\kern-.15em\lower.25ex\hbox{\the\scriptfont0 #2}}
\let\leukfrac\nfrac
\newcommand\half{\leukfrac1/2}
\newtheorem{theorem}{Theorem}[section]
\newtheorem{corollary}[theorem]{Corollary}
\newtheorem{lemma}[theorem]{Lemma}
\newtheorem{proposition}[theorem]{Proposition}
\numberwithin{equation}{section}
\theoremstyle{remark}
\newtheorem{remark}[theorem]{Remark}
\usepackage{epsfig}
\begin{document}

\title[Hereditary indecomposability]%
      {Remarks on hereditarily indecomposable continua}

\author{Klaas Pieter Hart}

\address{Faculty of Information Technology and Systems\\TU Delft\\
         Postbus 5031\\2600~GA {} Delft\\the Netherlands}
\email{k.p.hart@its.tudelft.nl}
\urladdr{http://aw.twi.tudelft.nl/\~{}hart}

\author{Jan van Mill}
\address{Faculty of Sciences\\
         Division of Mathematics\\
         Vrije Universiteit\\
         De Boelelaan 1081\textsuperscript{a}\\
         1081 HV Amsterdam\\
         The Netherlands}
\email{j.van.mill@cs.vu.nl}

\author{Roman Pol}
\address{Institute of Mathematics, University of Warsaw\\
         Banacha 2\\
         02-197 Warszawa\\
         Poland}
\email{pol@mimuw.edu.pl}

\date{\filedate}

\begin{abstract}
We recall a characterization of hereditary indecomposability originally
obtained by Krasinkiewicz and Minc, and show how it may be used to give
unified constructions of various hereditarily indecomposable continua.
In particular we answer a question asked by Mackowiak and Tymchatyn 
by showing that any continuum of arbitrary weight is a weakly confluent 
image of a hereditarily indecomposable continuum of the same weight. 

We present two methods of constructing these preimages:
(a)~by model-theoretic means, using the compactness and completeness
theorems from first-order logic to derive these results for continua of
uncountable weight from their metric counterparts; and
(b)~by constructing essential mappings from hereditarily indecomposable 
continua onto Tychonoff cubes. 

We finish by reviving an argument due to Kelley about hyperspaces
of hereditarily indecomposable continua and show how it leads to
a point-set argument that reduces Brouwer's Fixed-point theorem
to its three-dimensional version.
\end{abstract}



\maketitle

\section{Preliminaries}

\subsection{Hereditary indecomposable spaces}

A continuum is \emph{decomposable} if it can be written as the union of two 
proper subcontinua; it is \emph{indecomposable} otherwise.
A \emph{hereditarily indecomposable} continuum is one in which every 
subcontinuum is indecomposable.
It is easily seen that this is equivalent to saying that whenever two
continua in the space meet one is contained in the other.

This latter statement makes sense for arbitrary compact Hausdorff spaces,
connected or not; we therefore extend this definition and call a compact 
Hausdorff space \emph{hereditarily indecomposable} if it satisfies the 
statement above: whenever two continua in the space meet one is contained in 
the other.
Thus, zero-dimensional spaces are hereditarily indecomposable too.

We shall mainly use a characterization of hereditary indecomposability
that can be gleaned from \cite[Theorem~3]{KrasinkiewiczMinc1977} and which
was made explicit in~\cite[Theorem~2]{OversteegenTymchatyn86}.
To formulate it we introduce some terminology.

Let $X$ be compact Hausdorff and let $C$ and $D$ be disjoint closed subsets
of~$X$; as in~\cite{KrasinkiewiczMinc1977} we say that $(X,C,D)$~is
\emph{crooked} between the neigbourhoods $U$ of~$C$ and $V$ of~$D$ if we can 
write $X=X_0\cup X_1\cup X_2$, where each~$X_i$ is closed and, moreover,
$C\subseteq X_0$, $X_0\cap X_1\subseteq V$, $X_0\cap X_2=\emptyset$,
$X_1\cap X_2\subseteq U$ and $D\subseteq X_2$.
We say $X$ is crooked between $C$ and $D$ if $(X,C,D)$~is crooked between
any pair of neighbourhoods of~$C$ and~$D$.

We can now state the characterization of hereditary indecomposability
that we will use.

\begin{theorem}[Krasinkiewicz and Minc]
A compact Hausdorff space is hereditarily indecomposable if and only if it is 
crooked between every pair of disjoint closed (nonempty) subsets.
\end{theorem}

This characterization can be translated into terms of closed sets only;
we simply put $F=X\setminus V$ and $G=X\setminus U$, and reformulate
some of the premises and the conclusions.
We get the following formulation.

\begin{theorem}\label{thm.zig-zag}
A compact Hausdorff space $X$ is hereditarily indecomposable if and only if 
whenever four closed sets $C$, $D$, $F$ and $G$ in~$X$ are given such that 
$C\cap D=C\cap G=F\cap D=\emptyset$ one can write $X$ as the union of three 
closed sets $X_0$, $X_1$ and~$X_2$ such that
$C\subseteq X_0$,
\ $D\subseteq X_2$,
\ $X_0\cap X_1\cap G=\emptyset$,
\ $X_0\cap X_2=\emptyset$, and
  $X_1\cap X_2\cap F=\emptyset$.
\end{theorem}

To avoid having to write down many formulas we call a quadruple
$(C,D,F,G)$ with $C\cap D=C\cap F=D\cap G=\emptyset$ a
\emph{pliable foursome} and we call a triple $(X_0,X_1,X_2)$
with $C\subseteq X_0$,
\ $D\subseteq X_2$,
\ $X_0\cap X_1\cap G=\emptyset$,
\ $X_0\cap X_2=\emptyset$, and
  $X_1\cap X_2\cap F=\emptyset$
a \emph{chicane} for $(C,D,F,G)$.
Thus, a compact Hausdorff space is hereditarily indecomposable if and only if
there is a chicane for every pliable foursome.

This characterization can be improved by taking a base $\calB$ for the closed
sets of the space~$X$ that is closed under finite intersections.
The space is hereditarily indecomposable if and only if there is a chicane
for every pliable foursome whose terms come from~$\calB$.

To prove the nontrivial implication let $(C,D,F,G)$ be a pliable foursome
and let $(O_C,O_D,O_F,O_G)$ be a swelling of it, i.e.,
every set~$O_P$ is an open around~$P$ and $O_P\cap O_Q=\emptyset$ if and only
if $P\cap Q=\emptyset$, where $P$ and $Q$ run through $C$, $D$, $F$ and~$G$
(see \cite[7.1.4]{Engelking89}).
Now compactness and the fact that $\calB$ is closed under finite intersections
guarantee that there are $C'$, $D'$, $F'$ and~$G'$ in~$\calB$ such that
$P\subseteq P'\subseteq O_P$ for $P=C$, $D$, $F$,~$G$.
Any chicane for $(C',D',F',G')$ is a chicane for~$(C,D,F,G)$.

\subsection{A crooked partition of the square}
\label{subsec.crooked}

Let $P$ be the closure of the union of the five open rectangles in~$\I^2$,
depicted in Figure~\ref{fig.crooked} below. 
The set $\I^2\setminus P$ is the union of the disjoint open sets~$M_0$ 
and~$M_1$ in the picture. 
Clearly, $\cl{M}_0\cap \cl{M}_1=\emptyset$.
Observe that $\{0\}\times\I\subseteq M_0$ and $\{1\}\times\I\subseteq M_1$.
It follows that $P$ is a partition between $\{0\}\times\I$ and 
$\{1\}\times\I$ in~$\I^2$.

\begin{figure}[htb]
\centerline{\epsfig{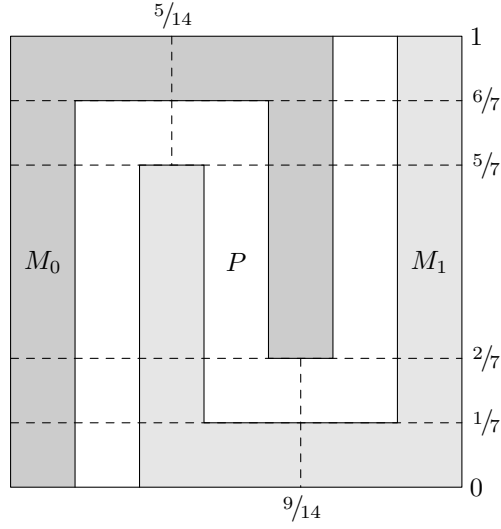}}
\caption{The crooked partition}
\label{fig.crooked}
\end{figure}

We shall use $P$ to create chicanes for pliable foursomes.
Here's how: given a pliable foursome $(C,D,F,G)$ apply Urysohn's lemma to get 
a continuous function $f:X\to\I$ such that
$f[C]=\{0\}$, $f[G]\subseteq[0,\half]$, $f[F]\subseteq[\half,1]$ and
$f[D]=\{1\}$.
One can then create a chicane by creating a continuous function $g:X\to\I$
such that $(g\mapdiagbin f)[X]\subseteq P$ and setting
$X_j=(g\mapdiagbin f)\preim[P_j]$, where 
$P_0=\bigl\{(x,y):x\le \leukfrac5/{14}\bigr\}$,
\ $P_1=\bigl\{(x,y):\leukfrac5/{14}\le x\le \leukfrac9/{14}\bigr\}$ and
$P_2=\bigl\{(x,y):\leukfrac9/{14}\le x\bigr\}$.
We shall call a function as $f$ a \emph{Urysohn function} for the 
foursome~$(C,D,F,G)$.

We summarize the foregoing discussion in the following lemmas.

\begin{lemma}\label{lemma.pliable.family}
Let $X$ be a compact space of weight~$\kappa$.
Then there is a family $\calF$ of continuous functions, from~$X$ to~$\I$,
of size~$\kappa$ such that for every pliable foursome $(C,D,F,G)$ there
is an~$f$ in~$\calF$ such that $f[C]=\{0\}$, $f[G]\subseteq[0,\half]$,
$f[F]\subseteq[\half,1]$ and $f[D]=\{0\}$. \qed
\end{lemma}

We call a family $\calF$ as in Lemma~\ref{lemma.pliable.family} a
\emph{pliable family} for~$X$.

\begin{lemma}\label{lemma.making.hi.spaces}
Let $X$ be a compact Hausdorff space and $\calF$ a pliable family of
functions for~$X$.
If $Y$ is a closed subspace of~$X$ with the property that for every~$f\in\calF$
there is a continuous function $g:X\to\I$ such that 
$(g\mapdiagbin f)[Y]\subseteq P$ then $Y$~is hereditarily indecomposable. 
\end{lemma}

\begin{proof}
Simply observe that every pliable foursome from~$Y$ is also a pliable
foursome in~$X$ and hence that the restrictions of the elements of~$\calF$
to~$Y$ form a pliable family for~$Y$.
\end{proof}

\subsection{Dimension and essential families}
\label{subsec.dim.and.essential}

We adopt the contrapositive of the Theorem on Partitions 
\cite[7.2.15]{Engelking89} as our definition of the covering dimension~$\dim$.
To this end we define a (finite or infinite) family 
$\bigl\{(A_i,B_i):i\in I\bigr\}$ of pairs of disjoint closed sets to be 
\emph{essential} if whenever we take partitions~$L_i$ between $A_i$ and~$B_i$
for all~$i$ the intersection~$\bigcap_{i\in I}L_i$ is nonempty.

If $X$ is a compact Hausdorff (or even normal) space and $n\in\N$ then we say
that $\dim X\ge n$ if $X$~has an essential family of pairs of closed sets
of cardinality~$n$; the \emph{covering dimension} of~$X$ is the maximum~$n$
such that $\dim X\ge n$, if such an~$n$ exists; we write $\dim X=\infty$
if $\dim X\ge n$ for all~$n$.
Note that $\dim X=\infty$ means that $X$ must have arbitrarily large
\emph{finite} essential families; if $X$ possesses an infinite
essential family then we say that $X$~is \emph{strongly infinite dimensional}.

The following lemma will be used to verify that certain spaces have a large
enough covering dimension.

\begin{lemma}\label{lemma.partial.essential}
Let $\bigl\{(A_i,B_i):i\in I\bigr\}$ be an essential family in a normal
space~$X$ and assume $I$ is split into two sets $J$ and~$K$.
Let, for $i\in J$, a partition~$L_i$ between $A_i$ and $B_i$ be given
and put~$L_J=\bigcap_{i\in J}L_i$.
Then $\bigl\{(A_i\cap L_J,B_i\cap L_J):i\in K\bigr\}$ is an essential family
in~$L_J$.
\end{lemma}

\begin{proof}
By normality we can extend, for every~$i\in K$, any partition in~$L_J$ 
between $A_i\cap L_J$ and $B_i\cap L_J$ to a partition in~$X$ between $A_i$ 
and~$B_i$.
Now apply the assumption that the full family is essential.
\end{proof}

\subsection{Faces of cubes}

Consider a Tychonoff cube~$\I^\kappa$.
For every~$\alpha\in\kappa$ we put
$A_\alpha=\{x:x_\alpha=0\}$ and $B_\alpha=\{x:x_\alpha=1\}$; these are
the $\alpha$th opposite faces of~$\I^\kappa$.
The following theorem is a fundamental fact about these faces.

\begin{theorem}\label{thm.tychonoff.essential}
The family $\bigl\{(A_\alpha,B_\alpha):\alpha\in\kappa\bigr\}$ is essential
in~$\I^\kappa$.
\end{theorem}

For finite $\kappa$ this follows from Brouwer's Fixed-Point Theorem
\cite[7.3.19]{Engelking89}.
In case $\kappa$~is infinite we put, for every finite subset~$a$ of~$\kappa$, 
$F_a=\pi_a\preim\bigl[\bigcap_{\alpha\in a}\pi_a[L_\alpha]\bigr]$,
where $\pi_a$ denotes the projection onto the subcube~$\I^a$.
By the finite case each~$F_a$ is nonempty and clearly $a\subseteq b$
implies $F_a\supseteq F_b$, so $\bigcap_aF_a\neq\emptyset$.
Now check that $\bigcap_a F_a=\bigcap_\alpha L_\alpha$.

With this fact in mind we call a continuous map $f:X\to\I^\kappa$ 
\emph{essential} if the family 
$\bigl\{(f\preim[A_\alpha],f\preim[B_\alpha]):\alpha<\kappa\bigr\}$
is essential.
A routine application of Urysohn's lemma shows that $X$~admits an essential
map onto~$\I^\kappa$ iff $X$ has an essential family of size~$\kappa$.

\subsection{Hyperspaces}

The \emph{hyperspace} of a space~$X$ is the family~$\HypX$ of nonempty
closed subsets of~$X$ endowed with the \emph{Vietoris topology}, which
has the family of sets of the form
$\Viet{U}=\{F:F\subseteq U\}$ and $\Viet{X,U}=\{F:F\cap U\neq\emptyset\}$,
where $U$~is open, as a subbase.
It is well-known that $\HypX$~is compact Hausdorff if $X$~is and that
if $X$~is compact metric, with metric~$d$, then the corresponding Hausdorff
metric~$d_H$ generates the Vietoris topology of~$\HypX$.
 
An important subspace of~$\HypX$ is $\CsX$, the space of all subcontinua
of~$X$; it is a closed subset of~$\HypX$, hence it is also compact
if $X$~is.

\section{Bing's continua}\label{bingeen}

We begin by constructing an infinite-dimensional hereditarily indecomposable
compact subset of the Hilbert cube~$\Hcube$.

To this end we let $\pi_i$ denote the projection of~$\Hcube$ onto the $i$-th
coordinate.
Furthermore we fix a pliable family $\{f_i:i\in\N\}$ of continuous functions 
for~$\Hcube$. 
For every $i$ we let $u_i=\pi_{2i}\mapdiagbin f_i$ be the diagonal map 
of~$\pi_{2i}$ and~$f_i$ from~$\Hcube$ to~$\I^2$.
The letter $P$ still refers to the partition of the square from 
Figure~\ref{fig.crooked}.

\begin{theorem}\label{thm.strongly.inf.dim.hi}
Let $X=\bigcap_{i=1}^\infty u_i\preim [P]$.
Then $X$ is an infinite-di\-men\-sional hereditarily indecomposable compact 
space.
\end{theorem}

\begin{proof}
Lemma~\ref{lemma.making.hi.spaces} implies immediately that $X$~is hereditarily
indecomposable: for every~$f_i$ the projection~$\pi_{2i}$  is as required.

To see that $X$ is infinite-dimensional we observe that $u_i\preim[P]$~is
a partition between the even-numbered faces~$A_{2i}$ and~$B_{2i}$ of~$\Hcube$ 
--- indeed:
$$
u_i\preim[P]\subseteq\{x\in\Hcube:\leukfrac1/7\le x_{2i}\le\leukfrac6/7\}.
$$
By Theorem~\ref{thm.tychonoff.essential} and 
Lemma~\ref{lemma.partial.essential} this implies that
$$
\bigl\{(A_{2i+1}\cap X,B_{2i+1}\cap X):i\in\N\}
$$
is an essential family in~$X$; so $X$~is even strongly infinite-dimensional.
\end{proof}

\begin{corollary}[Bing]\label{cor.Bing.n-dim}
For every $n$ there is an $n$-dimensional hereditarily indecomposable 
continuum.
\end{corollary}

\begin{proof}
As observed in the previous proof the traces of the odd-numbered faces 
of~$\Hcube$ on~$X$ form an essential family in~$X$.
One can therefore find a component~$S$ of~$X$ such that the traces from that 
family on~$S$ also form an essential family. 
Now let $\pi:S\to \I^{2n}$ be (the restriction of) the projection onto the 
first $2n$~coordinates.
Consider the monotone-light factorization of~$\pi$, i.e.,
write $\pi=\lambda\circ\mu$, where $\mu:S\to T$ is a monotone surjection 
and $\lambda:T\to \I^{2n}$ is a light map, cf.~\cite[6.2.22]{Engelking89}. 
Since $\lambda$~is light we have $\dim T\le 2n$,
cf.~\cite[7.4.20]{Engelking89}.

For odd $i<2n$ let $C_i=\lambda\preim[A_i]$ and $D_i=\lambda\preim[B_i]$ and 
observe that $\mu\preim[C_i]=S\cap A_i$ and $\mu\preim[D_i]=S\cap B_i$.
From these last equalities it follows that the~$C_i$ and~$D_i$ form an 
essential family in~$T$ and so~$\dim T\ge n$.

Because $n\le\dim T\le2n$ we may conclude that $T$~contains an $n$-dimensional
continuum~$B_n$.
Since $\mu$~is monotone and $S$~is hereditarily indecomposable, so is~$B_n$.
\end{proof}

\section{Bing's partitions}

We shall present a variation of the construction from the previous
section to demonstrate the following well-known result.

\begin{theorem}[Bing]\label{thm.hi.partition}
Let $X$ be a continuum and let $F_0$ and~$F_1$ be disjoint closed sets in~$X$.
Then there are disjoint open neighborhoods~$W_0$ and~$W_1$ of~$F_0$ and~$F_1$,
respectively, such that $X\setminus (W_0\cup W_1)$ is hereditarily 
indecomposable.
\end{theorem}

\begin{proof}
We use the partition~$P$ and the open sets~$M_0$ and~$M_1$ from
Section~\ref{subsec.crooked} again. 
Let $\{f_i:i\in\N\}$ be a pliable family for~$X$.

Choose open neighbourhoods $W_{0,0}$ and $W_{1,0}$ of~$F_0$ and $F_1$
respectively with disjoint closures.
Whenever $i\ge0$ and the open sets $W_{0,i}$ and $W_{1,i}$ with disjoint 
closures are found apply Urysohn's lemma to get a continuous function~$g_i$ 
such that $g_i[W_{0,i}]=\{0\}$ and $g_i[W_{1,i}]=\{1\}$ and set 
$$
    W_{0,i+1} = u_i\preim [M_0]\text{ and } W_{1,i+1} = u_i\preim [M_1],
$$
where $u_i=g_i\mapdiagbin f_i$.
Because the closures of $M_0$ and $M_1$ are disjoint the closures
of~$W_{0,i+1}$ and $W_{1,i+1}$ are disjoint as well.
Furthermore, because 
$u_i[W_{j,i}] \subseteq \{j\}\times \I\subseteq M_j$ 
we have $\cl{W}_{j,i}\subseteq W_{j,i+1}$ for $j=0$,~$1$.
In the end the sets 
$W_0 = \bigcup_{i=0}^\infty W_{0,i}$ and 
$W_1 = \bigcup_{i=0}^\infty W_{1,i}$
are disjoint open neighborhoods of $F_0$ and $F_1$, respectively. 

A direct application of Lemma~\ref{lemma.making.hi.spaces} shows that 
$L=X\setminus(W_0\cup W_1)$ is hereditarily indecomposable:
for every~$i$ the function~$g_i$ is a suitable partner for~$f_i$.
\end{proof}

\section{Continua of arbitrary weight}

This section contains some results on nonmetric continua.

\subsection{Bing's continua}

We begin by showing that nonmetric hereditarily indecomposable continua
of any prescribed weight exist.

\begin{theorem}\label{thm.unctble.existence.roman}
For every infinite~$\kappa$ there is a hereditarily indecomposable
continuum of weight~$\kappa$.
This continuum admits an essential map onto~$\I^\kappa$.
\end{theorem}

\begin{proof}
We use the proof of Theorem~\ref{thm.strongly.inf.dim.hi}.
Let $\{f_\alpha:\alpha\in\kappa\}$ be a pliable family of continuous
functions on the Tychonoff cube~$\I^\kappa$.
As before let $X=\bigcap_{\alpha\in\kappa}u_\alpha\preim[P]$,
where $u_\alpha=p_{2\alpha}\mapdiagbin f_\alpha$ and $p_\alpha$
is the projection onto the $\alpha$th coordinate.
The odd-numbered faces of~$\I^\kappa$ induce an essential family on~$X$;
it is also essential on some component of~$X$.
This component is the required continuum.
\end{proof}

\subsection{Bing's partitions}

There is no general analogue of Theorem~\ref{thm.hi.partition}; this
follows from the following (well-known) fact.

\begin{proposition}
Let $\L$ denote the long segment and let $L$ be any partition 
between $\{0\}\times\L$ and $\{1\}\times\L$ in the product $\I\times\L$;
then $L$ contains a copy of~$\I$.
\end{proposition}

\begin{proof}
Let $U$ and $V$ be disjoint open sets around $\{0\}\times\L$ and
$\{1\}\times\L$ respectively such that $L=(\I\times\L)\setminus(U\cup V)$.
Let $x_0=\sup\{x:(x,\omega_1)\in U\}$ and 
$x_1=\inf\{x:x>x_0$ and $(x,\omega_1)\in U\}$.
It is clear that $0<x_0\le x_1<1$ and that 
$[x_0,x_1]\times\{\omega_1\}\subseteq L$, so that we are done in case
$x_0<x_1$.
If $x_0=x_1$ then we can find an increasing sequence $\langle u_n\rangle_n$
and a decreasing sequence $\langle v_n\rangle_n$ such that
$(u_n,\omega_1)\in U$ and $(v_n,\omega_1)\in V$ for all~$n$.
Because the sets $U$ and $V$ are open and because $\omega_1$ has uncountable
cofinality we can find an~$\alpha<\omega_1$ such that 
$\{u_n\}\times[\alpha,\omega_1]\subseteq U$ and
$\{v_n\}\times[\alpha,\omega_1]\subseteq V$ for all~$n$.
It follows that 
$\{x_0\}\times[\alpha,\omega_1]\subseteq \cl U\cap\cl V\subseteq L$.
\end{proof}

Now consider $\I^{\omega_1}$ and embed $\L$ into $\I^{[1,\omega_1)}$;
this induces an embedding of $\I\times\L$ into~$\I^{\omega_1}$ so
that $\{0\}\times\L$ is embedded in the face~$A_0$ and $\{1\}\times\L$ is 
embedded in the face~$B_0$.
We see that every partition between~$A_0$ and~$B_0$ contains a copy of~$\I$.

\begin{remark}
An easy modification of the proof of Theorem~\ref{thm.hi.partition} will
show that in compact $F$-spaces of weight~$\aleph_1$ there are hereditarily
indecomposable partitions between any two disjoint closed sets.
Under the Continuum Hypothesis this applies to many \v{C}ech-Stone 
remainders such as $\beta\R^n\setminus\R^n$ and also
$\beta(\omega\times\I^{\omega_1})\setminus(\omega\times\I^{\omega_1})$.
\end{remark}

\subsection{Hereditarily indecomposable preimages}

In~\cite[(19.3)]{MackowiakTymchatyn1984} it is proven that every
metric continuum is the weakly confluent image of some hereditarily 
indecomposable metric curve.
A map is \emph{weakly confluent} if every continuum in the range is the
image of a continuum in the domain.

We shall show that this results holds in the nonmetric situation as well.

\begin{theorem}\label{thm.weakly.confluent.image.of.HI}
Every continuum is the continuous image of a one-dimensional hereditarily
indecomposable continuum (of the same weight) by a weakly confluent map.
\end{theorem}

For clarity of exposition we prove this theorem in stages;
first we show that every continuum is the continuous image of some
hereditarily indecomposable continuum of the same weight,
then we modify the construction to get a weakly confluent map and
finally we show how to make the domain one-dimensional.

\begin{proof}[Creating a hereditarily indecomposable preimage]
Let $X$ be a continuum of weight~$\kappa$ and assume that $X$~is embedded into
the Tychonoff cube~$\I^\kappa$.
Let $K$ be the hereditarily indecomposable continuum from 
Theorem~\ref{thm.unctble.existence.roman} and let $f:K\to\I^\kappa$ be an 
essential map.

For every finite subset $a$ of $\kappa$ we consider the map 
$\pi_a\circ f:K\to\I^a$ and the continuum~$\pi_a[X]$, where
$\pi_a$~is the projection of $\I^\kappa$ onto~$\I^a$.
Because $\pi_a\circ f$~is essential we may apply Theorem~4.3 
from~\cite{GrispolakisTymchatyn1981} to find a subcontinuum~$Y_a$
of~$K$ such that $(\pi_a\circ f)[Y_a]=\pi_a[X]$.
Because $K$~is compact the net $\langle Y_a:a\in[\kappa]^{<\omega}\rangle$
has a convergent subnet in~$\Cs{K}$; its limit~$Y$ is a subcontinuum of~$K$
that maps onto~$X$.
\end{proof}

To be able to improve this construction so as to make the map weakly
confluent we need the following result, which follows from
Theorem~3.5 of~\cite{GrispolakisTymchatyn1978}.

\begin{lemma}\label{lemma.G-P.3.5}
Let $n\in\N$ and let $X$ be a continuum in $\I^n\times\{0\}$ then there
is a copy~$H$ of the half~line~$[0,\infty)$ in~$\I^n\times(0,1]$ such that
$H\cup X=\cl H$ is a continuum with the property that for every
continuous surjection $f:Z\to \cl H$, where $Z$~is a continuum,
there is a subcontinuum~$Y$ of~$Z$ such that $f\restr Y:Y\to X$ is onto
and weakly confluent.
\end{lemma}

Using this lemma and a slightly more complicated proof we can ensure
that we get a weakly confluent map from a hereditarily indecomposable
continuum onto a given one.

\begin{proof}[Creating a weakly confluent preimage]
We now assume that our continuum~$X$ is embedded in~$\I^\kappa\times\{0\}$,
which we consider to be a subset of~$\I^\kappa\times\I$.
We take the continuum~$K$ from Theorem~\ref{thm.unctble.existence.roman}
and let $f:K\to\I^\kappa\times\I$ be an essential map.
For every finite subset~$a$ of~$\kappa$ we let $\pi_a$ denote the projection
of $\I^\kappa\times\I$ onto~$\I^a\times\I$.

An application of Lemma~\ref{lemma.G-P.3.5} yields for every finite set~$a$
a subcontinuum~$Y_a$ of~$K$ such that $(\pi_a\circ f)[Y_a]=\pi_a[X]$ and
the restriction $(\pi_a\circ f)\restr Y_a$ is weakly confluent.
As before we take a convergent subnet $\{Y_\alpha:\alpha\in A\}$ 
of $\{Y_a:a\in[\kappa]^{<\omega}\}$ with limit~$Y$; then $f[Y]=X$ and it 
remains to show that $f\restr Y$ is weakly confluent.
To this end let $C$ be a subcontinuum of~$X$ and choose for every~$a$ a
subcontinuum~$D_a$ of~$Y_a$ such that $(\pi_a\circ f)[D_a]=\pi_a[C]$.
The subnet $\{D_\alpha:\alpha\in A\}$ of~$\{D_a:a\in[\kappa]^{<\omega}\}$
has a convergent subnet $\{D_\beta:\beta\in B\}$ with limit~$D$;
it should be clear that $D\subseteq Y$ and $f[D]=C$.
\end{proof}

Finally we indicate how to get a \emph{one-dimensional} hereditarily
indecomposable continuum~$Y_1$ and a weakly confluent map from~$Y_1$
onto~$X$.
For this we need the following lemma.

\begin{lemma}
There are a one-dimensional subcontinuum~$U$ of~$\I^\kappa$ and a monotone
surjection $m:U\to\I^\kappa$.
\end{lemma}

\begin{proof}
This is a straightforward generalization of the proof of~19.1 
in~\cite{MackowiakTymchatyn1984}.
Let $\Cantor$ denote the standard Cantor~set in~$\I$.
For $\alpha\in\kappa$ put $U_\alpha=\{x\in\I^\kappa:{}$if $\beta\neq\alpha$
then~$x_\beta\in\Cantor\}$ and set $U=\bigcup_{\alpha\in\kappa}U_\alpha$.
Just as in~\cite{MackowiakTymchatyn1984} one verifies that $U$~is a closed
and connected subset of~$\I^\kappa$; to see that $U$~is one-dimensional
one only has to realize that every basic open cover lives on a finite
subset of~$\kappa$ and hence that it can be given an open refinement
of order~$2$.
Finally, the map $h^\kappa:U\to\I^\kappa$, where $h:\I\to\I$ is the
Cantor step~function, is a monotone map.
\end{proof}

\begin{proof}[A one-dimensional preimage]
By the previous lemma we can find a one-dimensional continuum~$X_1$
of the same weight as~$X$ and a monotone surjection $m:X_1\to X$.
Next find a hereditarily indecomposable continuum~$Y$ and a weakly confluent
surjection $f:Y\to X_1$.
As in the proof of Corollary~\ref{cor.Bing.n-dim} we take the monotone-light 
factorization of~$f$, i.e., a space~$Y_1$, a monotone map $\mu:Y\to Y_1$ and 
a light map $\lambda:Y_1\to X_1$ so that $f=\lambda\circ\mu$.
Because $\mu$~is monotone the space~$Y_1$ is hereditarily indecomposable,
because $\lambda$~is light it is one-dimensional and because $f$~is
weakly confluent so is~$\lambda$ and because $m$~is monotone the map 
$m\circ\lambda$ is weakly confluent.
\end{proof}

\subsection{Hereditarily infinite-dimensional spaces}

There are many constructions of hereditarily infinite-dimensional
continua, i.e., continua with only in\-fi\-nite-di\-men\-sional 
(nontrivial) subcontinua.
One such construction, due to Levin~\cite{Levin1995}, 
admits a striking generalization to higher cardinals.

\begin{theorem}
For every cardinal~$\kappa$ there is a hereditarily indecomposable
space of weight~$\kappa$ such that every subcontinuum of it has an essential
family of cardinality~$\kappa$.
\end{theorem}

\begin{proof}
Consider the continuum~$K$ constructed in 
Theorem~\ref{thm.unctble.existence.roman} and its essential family 
$E=\bigl\{(A_\alpha,B_\alpha):\alpha\le\kappa\bigr\}$.
Partition $\kappa$ into $\kappa$~many sets~$T_\alpha$ of size~$\kappa$
and let $\{B_\alpha:\alpha<\kappa\}$ be a base for~$K$.

For each~$\alpha$ let $W_\alpha$ be the union of all components 
of~$\cl B_\alpha$ on which 
$E_\alpha=\bigl\{(A_\beta,B_\beta):\beta\in T_\alpha\bigr\}$
is not essential.
Observe that $W_\alpha$~is open in~$\cl B_\alpha$ and that the 
family~$E_\alpha$ is not essential on any compact subset of~$W_\alpha$.

Next let $U_\alpha=W_\alpha\cap B_\alpha$ for each $\alpha$ and put
$U=\bigcup_{\alpha<\kappa}U_\alpha$.
The set~$U$ is open in~$K$ and the family~$E$ is not essential on any
compact subset of~$U$: if $C$~is such a set cover it by finitely 
many~$U_\alpha$ and use the disjointness of the sets~$T_\alpha$ to make
a set of partitions whose intersection misses~$C$.

It follows that every partition between $A_\kappa$ and $B_\kappa$ must meet
$K\setminus U$ and hence that $K\setminus U$~contains a non-trivial 
continuum~$H$.
Let $C$~be any subcontinuum of~$H$, let $p\in C$ and fix~$\alpha$
such that $p\in B_\alpha$ and $C\nsubseteq \cl B_\alpha$.
Consider the component~$Z$ of~$p$ in~$\cl B_\alpha$; because $K$~is
hereditarily indecomposable we have $Z\subseteq C$.
But then $Z\cap W_\alpha=\emptyset$ and so $E_\alpha$~is essential on~$Z$
and hence on~$C$.
\end{proof}

\section{Model-theoretic considerations}

In this section we call attention to the curious fact that many results
about compact spaces of uncountable weight can be derived by model-theoretic
means; in fact, the Compactness Theorem and the L\"owenheim-Skolem Theorem
enable one to deduce the uncountable versions directly from the theorems in 
the metric case.

\subsection{Wallman spaces}

The basis for the model-theoretic approach is Wallman's generalization
\cite{Wallman38} of Stone's representation theorem for Boolean algebras
to distributive lattices.
If $L$ is a distributive lattice (with $\0$ and $\1$) then there is
is a compact $T_1$-space~$wL$ with a base for its closed sets that is a 
homomorphic image of\/~$L$.
The homomorphism is an isomorphism if and only if $L$ is \emph{disjunctive},
which means: if $a\notle b$ then there is $c\in L$ such that $c\le a$
and $c\meet b=\0$.
Every compact $T_1$-space~$X$ can be obtained in this way: $X$~is the Wallman
space of its own family of closed sets.
From this it is clear that $wL$~is not automatically Hausdorff; in fact
$wL$~is Hausdorff if and only if $L$~is \emph{normal}, which is expressed
as follows:
\begin{equation}\label{L.is.normal}
(\forall x)(\forall y)(\exists u)(\exists v)
\bigl[(x\meet y=\0)\implies
\bigl((x\meet u=\0)\land(y\meet v=\0)\land(u\join v=\1)\bigr)\bigr].
\end{equation}
In a similar fashion we can express that $wL$~is connected or hereditarily
indecomposable.
The following formula expressed the connectivity of~$wL$:
\begin{equation}\label{L.is.connected}
(\forall x)(\forall y)
\bigl[\bigl((x\meet y=\0)\land(x\join y=\1)\bigr)\implies
      \bigl((x=\0)\lor(x=\1)\bigr)\bigr].
\end{equation}
This suffices because every base for the closed sets of a compact space
that is a lattice contains every clopen set of that space.
We can read this formula as expressing a property of~$\1$, to wit
``$\1$ is connected''; we therefore abbreviate it as~$\conn(\1)$
and we shall write $\conn(a)$ to denote Formula~\ref{L.is.connected}
with $\1$ replaced by~$a$ and use it to express that $a$~is connected
(or better: the set represented by~$a$ is connected).

To ensure that $wL$~is hereditarily indecomposable it suffices to have
a chicane for every pliable foursome from~$L$ and this is exactly
what the following formula expresses.
\begin{multline}\label{L.is.HI}
(\forall x)(\forall y)(\forall u)(\forall v)(\exists z_1,z_2,z_3)
\bigl[\bigl((x\meet y=\0)\land(x\meet u=\0)\land(y\meet v=\0)\bigr)\implies{}\\
{}\implies\bigl((x\meet(z_2\join z_3)=\0)\land (y\meet(z_1\join z_2)=\0)
   \land (z_1\meet z_3=\0){}\\
{}\land(z_1\meet z_2\meet v=\0)\land(z_2\meet z_3\meet u=\0)\land
 (z_1\join z_2\join z_3=\1)      \bigr)\bigr].
\end{multline}

\subsection{Existence of hereditarily indecomposable continua}

The existence of the pseudoarc~$\psarc$ implies that there are one-dimensional
hereditarily indecomposable continua of arbitrarily large weight.
Indeed, the family of closed sets of~$\psarc$ is a distributive and disjunctive
lattice that satisfies formulas~\ref{L.is.normal}, \ref{L.is.connected}
and~\ref{L.is.HI}; it also satisfies
\begin{multline}\label{dimL<=1}
(\forall x_0)(\forall y_0)(\forall x_1)(\forall y_1)
(\exists u_0,v_0,u_1,v_1)
\bigl[\bigl((x_0\meet y_0=\0)\land(x_1\meet y_1=\0)\implies{}\\
{}\implies\bigl((x_0\meet u_0=\0)\land (y_0\meet v_0=\0)\land
                (x_1\meet u_1=\0)\land (y_1\meet v_1=\0)\land{}\\
        {}\land (u_0\join v_0=\1)\land (u_1\join v_1=\1)\land
                (u_0\meet v_0\meet u_1\meet v_1=\0)\bigr)\bigr].
\end{multline}
This formula expresses $\dim wL\le1$ in terms of closed sets, see
Section~\ref{subsec.dim.and.essential}.
Therefore this combination of formulas is consistent and so, by the (upward)
L\"owenheim-Skolem theorem, it has models of every cardinality.
Thus, given a cardinal~$\kappa$ there is a distributive and disjunctive
lattice~$L$ of cardinality~$\kappa$ that satisfies~\ref{L.is.normal}, 
\ref{L.is.connected}, \ref{L.is.HI} and~\ref{dimL<=1}.
The space~$wL$ is compact Hausdorff, connected, hereditarily indecomposable,
one-dimensional and of weight~$\kappa$ or less, but with at least~$\kappa$
closed sets.
Thus, if $\kappa\ge2^\lambda$ then the weight of~$wL$ is at least~$\lambda$.

To get a space of weight exactly~$\kappa$ we make sure that $wL$ has at least
$2^\kappa$ many closed sets.
To this end we introduce two sets of $\kappa$~many constants
$\{a_\alpha:\alpha<\kappa\}$ and $\{b_\alpha:\alpha<\kappa\}$ and two sets
of $\kappa$~many formulas:
for every~$\alpha$ the formula $a_\alpha\meet b_\alpha=\0$ and for
any pair of disjoint finite subsets~$p$ and~$q$ of~$\kappa$ the formula
$\bigmeet_{\alpha\in p}a_\alpha\meet\bigmeet_{\alpha\in q}b_\alpha\neq\0$.
Thus we have expanded the language of lattices by a number of constants
and we have added a set of formulas to the formulas that we used above.
This larger set~$\calT_\kappa$ of formulas is still consistent.

Take a finite subset~$T$ of $\calT_\kappa$ and fix a finite subset~$t$ 
of~$\kappa$ such that whenever $a_\alpha\meet b_\alpha=\0$ or
$\bigmeet_{\alpha\in p}a_\alpha\meet\bigmeet_{\alpha\in q}b_\alpha\neq\0$
belong to~$T$ we have $\alpha\in t$ and $p\cup q\subseteq t$.
Now take a map~$f$ from~$\psarc$ onto the cube~$\I^t$ and interpret
$a_\alpha$ by~$f\preim[A_\alpha]$ and $b_\alpha$ by~$f\preim[B_\alpha]$;
in this way we have ensured that every formula from~$T$ holds in the family
of closed subsets of~$\psarc$.
Therefore $T$~is a consistent set of formulas and so, because it was arbitrary
and by the compactness theorem, the full set~$\calT_\kappa$ is 
consistent.

Because $\calT_\kappa$ has cardinality~$\kappa$ it has a model~$L$
of cardinality~$\kappa$.
Now $wL$~is as required: its weight is at most~$\kappa$ because $L$~is a base
of cardinality~$\kappa$.
On the other hand: for every subset~$S$ of~$\kappa$, we have, by compactness,
a nonempty closed set
$$
F_S=\bigcap_{\alpha\in S}a_\alpha\cap\bigcap_{\alpha\notin S}b_\alpha
$$
such that $F_S\cap F_T=\emptyset$ whenever $S\neq T$.

\begin{remark}
The reader may enjoy modyfying the above argument so as to ensure
that $\bigl\{(a_\alpha,b_\alpha):\alpha<\kappa\bigr\}$ is an essential family
in~$wL$.
To this end write down, for every finite subset~$a$ of~$\kappa$,
a formula~$\phi_a$ that expresses that 
$\bigl\{(a_\alpha,b_\alpha):\alpha\in a\bigr\}$
is essential.
Theorem~\ref{thm.strongly.inf.dim.hi} more than ensures that the
set of formulas consisting of~\ref{L.is.normal}, \ref{L.is.connected}, 
\ref{L.is.HI} and the~$\phi_a$ is consistent. 
\end{remark}

\subsection{Hereditarily indecomposable preimages}

We can also give a model-the\-o\-retic proof of 
Theorem~\ref{thm.weakly.confluent.image.of.HI}.
For this we need to know how to ensure that $wL$ maps onto the given
continuum and how to get this surjection to be weakly confluent.

\subsubsection*{Making a continuous surjection}

The following lemma tells us how to make continuous surjections.

\begin{lemma}
Let $X$ be compact Hausdorff and $L$ some normal, distributive and disjunctive
lattice.
If $X$ has a base~$\calB$ for the closed sets that is a lattice and embeddable
into~$L$ then $wL$ admits a continuous surjection onto~$X$.
\end{lemma}

\begin{proof}
We only sketch the argument.
Let $\phi:\calB\to L$ be an embedding and define $f:wL\to X$ by
``$f(p)$~is the unique point in 
$\bigcap\bigl\{C\in\calB:p\in\phi(C)\bigr\}$''.
It is straightforward to check that $f$~is onto and that
$f\preim[C]=\phi(C)$ for all~$C$.
\end{proof}

This tells us that to get a (one-dimensional) hereditarily indecomposable
continuum that maps onto the given continuum~$X$ we need to construct a 
distributive, disjunctive and normal lattice~$L$ that satisfies 
formulas~\ref{L.is.connected} and~\ref{L.is.HI} (and~\ref{dimL<=1}), 
and an embedding~$\phi$ of some base~$\calB$ for the closed sets of~$X$ 
into~$L$.

Let a continuum~$X$ and a lattice-base~$\calB$ for its closed sets be given.
As before we start with the formulas that ensure that $wL$~will be a 
hereditarily indecomposable continuum.
To these formulas we add the diagram of~$\calB$; this consists of~$\calB$
itself, as a set of constants, and the `multiplication tables' for~$\meet$
and~$\join$, i.e., $A\meet B=C$ whenever $A\cap B=C$ and $A\join B=C$
whenever $A\cup B=C$.

Now, if $L$~is to satisfy the diagram of~$\calB$ it must contain 
elements~$x_A$ for every $A\in\calB$ so that $x_{A\meet B}=x_A\meet x_B$
and $x_{A\join B}=x_A\join x_B$ hold whenever appropriate; but this simply
says that there is an embedding of $\calB$ into~$L$.

We are left with the task of showing that the set~$\calT$ of formulas that 
express distributivity, disjunctiveness, normality as well as
formulas~\ref{L.is.connected} and~\ref{L.is.HI} (and~\ref{dimL<=1}), 
together with the diagram of~$\calB$ is consistent.
Let $T$ be a finite subset of~$\calT$ and, if necessary, add the first six 
formulas to it.
Let $\calB'$ be a countable, normal and disjunctive sublattice of~$\calB$
that contains the finitely many constants that occur in~$T$.
The Wallman space of~$\calB'$, call it~$Y$, is a metric continuum
and therefore the continuous image of a hereditarily indecomposable
(one-dimensional) continuum~$K$.
The lattice of closed sets of~$K$ satisfies all the formulas from~$T$:
interpret $A$ by its preimage in~$K$.

It follows that $\calT$~is consistent and that it therefore has a model~$L$
of the same cardinality as~$\calT$, which is the same as the cardinality
of~$\calB$.
The lattice~$L$ satisfies all formulas from~$\calT$; its Wallman space
is a (one-dimensional) hereditarily indecomposable continuum that maps
onto~$X$.
If $\calB$ is chosen to be of minimal size then $wL$ is of the same weight
as~$X$.

This proof is much like the model-theoretic proof of the main theorem
of~\cite[Section~2]{DowHart2000b} which says that every continuum
of weight~$\aleph_1$ is a continuous image of the \v{C}ech-Stone
remainder of the real line.

\subsubsection*{Making a weakly confluent map}

We now improve the foregoing construction so as to make the continuous
surjection weakly confluent.

The following theorem --- which is a souped-up version of the Marde\v{s}i\'c
factorization theorem --- implies that it suffices to get some hereditarily 
indecomposable continuum~$Y$ that admits a weakly confluent map~$f$ 
onto our continuum~$X$.

\begin{theorem}
Let $f:Y\to X$ be a continuous surjection between compact Hausdorff spaces. 
Then $f$ can be factored as~$h\circ g$, where 
$Y\buildrel g\over\to Z\buildrel h\over\to X$ and $Z$ has the same weight 
as~$X$ and shares many properties with~$Y$.
\end{theorem}

\begin{proof}
Let $\calB$ be a lattice-base for the closed sets of~$X$ (of minimal size) 
and identify it with its copy $\{f\preim[B]:B\in\calB\}$ in~$\Hyp{Y}$.
By the L\"owenheim-Skolem theorem 
\cite[Corollary 3.1.5]{Hodges93} there is a lattice~$\calD$, of the same
cardinality as~$\calB$, such that $\calB\subseteq\calD\subseteq\Hyp{Y}$ 
and $\calD$~is an elementary substructure of~$\Hyp{Y}$.
The space~$Z=w\calD$ is as required.
\end{proof}

Some comments on this theorem and its proof are in order, because they do
not seem to say very much.
However, `elementary substructure' is an extremely powerful concept.
In our context it means that the smaller structure is closed off under every 
possible finitary lattice-theoretic operation of interest.

For example, if $Y$ hereditarily indecomposable then $\Hyp{Y}$ is closed
under the operation, implicit in formula~\ref{L.is.HI}, which assigns
a chicane to every pliable foursome.
But then $\calD$ must be closed under this operation as well and hence
$Z$~is hereditarily indecomposable.

Likewise $\dim Z=\dim Y$, because if there is an essential family in~$\Hyp{Y}$
of size~$n$ then there must be one in~$\calD$ (use a constant operation that
assigns an essential family of size~$n$ to everything) and, conversely, 
if there is an essential family of size~$n$ in~$\calD$ then it is essential
in~$\Hyp{Y}$ as well: $\calD$~is closed under the operation of assigning
sequences of partitions with empty intersection to inessential families.

We leave to the reader the verification that if $f$ is weakly confluent
then so is the map~$h$ in the factorization.

Now let $X$ be a continuum.
Our aim is of course to find a lattice~$L$ that contains the diagram 
of~$\HypX$ --- to get our continuous surjection~$f$ --- 
and for every~$C\in\CsX$ a continuum~$C'$ in~$wL$ such that $f[C']=C$.

As before we add the diagram of~$\HypX$ to the formulas that guarantee
that $wL$ will be a hereditarily indecomposable continuum.
In addition we take a set of constants $\{C':C\in\CsX\}$ and stipulate
that $C'$ will be a continuum that gets mapped onto~$C$.

To make sure that every $C'$ is connected we put $\conn(C')$ into
our set of formulas, for every~$C$.
Next, $f[C']\subseteq C$ translates, via the embedding into~$L$, into
$C'\le C$ (or better $C'=C'\meet C$).
Now, if it happens that $f[C']\subsetneq C$ then there is a closed set~$D$
in~$X$ (in fact it is $f[C']$ but that is immaterial) such that
$C'\le D$ and $C\notle D$.
In order to avoid this we also add, for every~$C\in\CsX$ and 
every~$D\in\HypX$, the formula
$$
(C'\le D) \implies (C\le D)
$$
to our set of formulas.

Again, the theorem in the metric case implies that this set of formulas is
consistent --- given a finite subset~$T$ of it make a metric 
continuum~$X_T$ as before, by expanding $\{B\in\HypX:B$~occurs in~$T\}$
to a countable normal sublattice~$\calB$ of~$\HypX$; 
then find a metric continuum~$Y_T$ of the desired type that admits a weakly
confluent map~$f$ onto~$X_T$; finally choose for every 
$C\in\CsX$ that occurs in~$T$ a continuum in~$Y_T$ that maps onto~$C$ and 
assign it to~$C'$; this then makes the family of closed sets of~$Y_T$ a model
of~$T$.

As before we obtain a lattice~$L$ whose Wallman space is one-di\-men\-sional
and hereditarily indecomposable, and which, in addition, admits a weakly confluent
map onto~$X$.

\section{From three to infinity}

In this section we shall show that Brouwer's Fixed-point Theorem in dimension
three implies all of its higher-dimension versions, using only point-set
arguments and a smattering of Linear Algebra.
The point-set arguments can be culled from Kelley's proof, 
from~\cite{Kelley42}, of his theorem the hyperspace of a (at least) 
two-dimensional hereditarily indecomposable continuum is infinite-dimensional.
To convince the reader that point-set arguments really suffice and to make
Kelley's result better known we shall give the argument in full.
In this section all continua under consideration are metrizable;
we invariably use $\rho$ to denote a compatible metric and
$\rho_H$ to denote the corresponding Hausdorff metric.

\subsection{More on hyperspaces}

Most of our arguments will take place in the hyperspace~$\CsX$ of all 
subcontinua of a two-dimensional hereditarily indecomposable continuum~$X$.

\subsubsection*{Order arcs}

It is well-known that $\CsX$ is arcwise connected whenever $X$~is a metric
continuum; in fact if $A\in\CsX$ then there is a linearly ordered 
family~$\calC$ of continua containing~$A$ and~$X$ and that is homeomorphic
to~$\I$.
For hereditarily indecomposable continua we can give a completely elementary
proof of this fact.

\begin{lemma}\label{lemma.Cx.is.arc}
Let $X$ be a hereditarily indecomposable continuum and for $x\in X$ put
$\calC_x=\{C\in\CsX:x\in C\}$.
Then $\calC_x$ is a chain, whose subspace and order topologies coincide and 
make it homeomorphic to~$\I$.
\end{lemma}

\begin{proof}
That $\calC_x$ is a chain follows from hereditary indecomposability of~$X$.
It is clear that $\calC_x$ is complete: if $\calF\subseteq\calC_x$ then
$\cl\bigcup\calF$ is the supremum of~$\calF$ in~$\calC_x$.
To see that $\calC_x$ has no jumps take $C$ and $D$ in~$\calC$ 
with $C\subsetneq D$ and fix an open set~$U$ such that $C\subseteq U$ and 
$D\nsubseteq\cl U$.
Now the component~$E$ of $\cl U$ that contains~$C$ meets the boundary of~$U$,
so $C\subsetneq E$, and is contained in~$\cl U$, so $E\subsetneq D$.

The set $\calC_x$~is closed in~$\CsX$: its complement,
$\bigl\{C:C\subseteq X\setminus\{x\}\bigr\}$, is a basic open set.
Likewise the sets $\{C:C\subseteq A\}$ and $\{C:A\subseteq C\}$ are closed
in~$\CsX$; this shows that the order topology on~$\calC_x$ is contained
in the subspace topology.
Because both topologies are compact Hausdorff they coincide; because
this topology is metric we find that $\calC_x$~is isomorphic and homeomorphic
to~$\I$.
\end{proof}

\subsubsection*{Whitney levels}

A Whitney map for $\HypX$ is a continuous function $\mu:\HypX\to\R$ such 
that $\mu\bigl(\{x\}\bigr)=0$ for all~$x$ and $\mu(C)<\mu(D)$ whenever 
$C\subsetneq D$.
If $X$~is compact metric then there are Whitney maps for~$\HypX$,
see \cite[4.33]{Nadler1992}.

We fix a hereditarily indecomposable continuum~$X$ and a Whitney map
$\mu:\CsX\to\R$ (we shall work inside $\CsX$ only).
The fibers $\mu\preim(r)$ ($0\le r\le\mu(X)$) divide $\CsX$ into layers,
refered to as Whitney levels.
We list some properties of Whitney levels.

\begin{lemma}
Every Whitney level is closed. \qed
\end{lemma}

\begin{lemma}
Every Whitney level is a pairwise disjoint family of continua.
\end{lemma}

\begin{proof}
Apply hereditary indecomposability.
\end{proof}

\begin{lemma}
Every Whitney level covers $X$.
\end{lemma}

\begin{proof}
The function $\mu$ is continuous and, for every~$x$, the set~$\calC_x$
is an arc that connects~$\{x\}$ and~$X$; it follows that
$\mu[\calC_x]=\bigl[0,\mu(X)\bigr]$.
\end{proof}

One can easily show that $\rho_H(A,B)<\delta$ implies
$\bigl|\diam(A)-\diam(B)\bigr|<2\delta$, so that $\diam$ is a continuous
function on~$\HypX$.
It follows that, for every~$r$, the diameter function assumes a minimum on
the Whitney level~$\wlev{r}$.
On the other hand, for every positive number~$\epsilon$ the 
set~$\{A:\diam(A)<\epsilon\}$ is an open neighbourhood of the closed set 
$\bigl\{\{x\}:x\in X\bigr\}$; it follows, by compactness, that there is 
a positive number~$s$ such that $\wlev{r}\subseteq\{A:\diam(A)<\epsilon\}$
whenever $r<s$.

We now have all the ingedients we need to be able to present Kelley's
argument.

\subsection{Kelley's argument}

For the remainder of this section we fix a hereditarily indcomposable
continuum~$X$ that is at least two-di\-men\-sional and we fix an essential
family $\bigl\{(A_0,B_0),(A_1,B_1)\bigr\}$ witnessing this.

To begin fix $\epsilon>0$ such that $\rho(x,y)>\epsilon$ whenever
$x\in A_i$ and $y\in B_i$, where $i=0$,~$1$.
We may assume, without of loss of generailty, that $\epsilon=1$
(if necessary scale $\rho$ by the factor~$\frac1\epsilon$).
The following lemma will be used toward the end of Kelley's argument.

\begin{lemma}\label{lemma.small.mesh}
Let $\calN$ be a finite disjoint collection of closed sets with 
diameter at most~$\half$ in~$X$. 
Then there is a continuum in~$X$ of diameter at least\/~$1$ that misses
all elements from~$\calN$.
\end{lemma}

\begin{proof}
Striving for a contradiction, assume that for $\calN$ such as in the 
formulation of the lemma, every continuum in~$X\setminus\bigcup\calN$ 
has diameter less than~$1$. 
Since $\calN$ is finite, there clearly is a finite disjoint collection 
$\calN_0$ of closed subsets of~$X$ of mesh less than~$1$ such that 
$\bigcup\calN$ is contained in the interior~$W$ of~$\bigcup\calN_0$. 
Since by assumption each component of~$X\setminus W$ has diameter less 
than~$1$, the set~$X\setminus W$ can also be covered by a finite 
disjoint collection~$\calN_1$ of closed sets with diameter less 
than~$1$. 
For $i=0$ and~$1$ let $S_i$ be a closed set in~$X$ separating the disjoint
closed sets
\begin{align*}
C_i&=A_i\cup\bigcup\{N\in\calN_i:N\cap B_i=\emptyset\}\rlap{ and}  \\
D_i&=B_i\cup\bigcup\{N\in\calN_i:N\cap B_i=\emptyset\}.
\end{align*}
Then, clearly, $S_0\cap S_1=\emptyset$, which in turn contradicts our 
assumption that the pairs $(A_0,B_0)$ and $(A_1,B_1)$ form an essential
family.
\end{proof}

Let us now take a Whitney map $\mu:\CsX\to\R$ and fix $s>0$ such that
the Whitney level $\wlev{r}$ is contained in $\{A:\diam(A)<\half\}$,
whenever $r<s$.
We shall show that $\wlev{r}$ is infinite-dimensional whenever $r<s$.
Fix such an~$r$ and put $\eta=\min\{\diam(A):\mu(A)=r\}$.
The following proposition implies that the Whitney level $\wlev{r}$
is infinite dimensional --- we shall explain this later in 
Remark~\ref{rem.wlev(r).inf.dm}.

\begin{proposition}\label{prop.many.neighbours}
Every finite closed cover of $\wlev{r}$ of mesh less than~$\eta/(4n)$
has an element that meets at least\/ $n$ other elements of the cover.
\end{proposition}

\begin{proof} 
Put $\epsilon=\eta/(4n)$ and assume that $\mathfrak{F}$ is a finite closed
cover of~$\wlev{r}$ with mesh less than~$\epsilon$ such that each element 
of~$\mathfrak{F}$ meets at most $n-1$ other elements of~$\mathfrak{F}$. 
We shall associate to each element $\calA\in\mathfrak{F}$ a compact 
subset~$\varphi(\calA)$ of~$\bigcup\calA$
such that
\begin{enumerate}
\item $\varphi(\calA)$ meets every element of~$\calA$,
\item $\diam\varphi(\calA)\le2\epsilon$,
\item $\varphi(\calA)\cap \varphi(\mathcal{B})=\emptyset$ whenever
      $\calA$ and $\calB$ are distinct elements of~$\mathfrak{F}$.
\end{enumerate}
Assume that $\varphi$ is already defined on a subfamily~$\mathfrak{G}$ 
of~$\mathfrak{F}$, and take $\calA$ in $\mathfrak{F}\setminus \mathfrak{G}$;
we show how to extend $\varphi$ to~$\mathfrak{G}\cup\{\calA\}$ (this then
means that we can define $\varphi$ on all of~$\mathfrak{F}$ in finitely
many steps).

The set $\mathfrak{H}$ of all elements of $\mathfrak{G}$ that meet~$\calA$ has,
by assumption, cardinality less than~$n$; because 
$\calA\cap\calB\neq\emptyset$ iff $\bigcup\calA\cap\bigcup\calB\neq\emptyset$
it suffices sure that $\varphi(\calA)$ does not meet~$\varphi(\calB)$
for any~$\calB$ in~$\mathfrak{H}$.
 
For each~$\calB$ in~$\mathfrak{H}$ let $B(\calB)$ be the closed $\epsilon$-ball
about $\varphi(\calB)$ and fix an $A\in \calA$. 
We shall show that $\{B(\calB):\calB\in\mathfrak{H}\}$ does not cover~$A$.
Indeed, otherwise, because $A$~is connected, we could arrange this family
into a sequence $B_1,\ldots,B_p$, $p< n$, with 
$B_i\cap \bigcup_{j<i} B_j\neq\emptyset$ for each $i\le p$. 
But then we could find an upper bound for $\diam A$, thus:
$$
\diam A\le (n-1)\cdot \max_{i\le p}\diam B_i
       \le (n-1)\cdot (4\epsilon) <\eta.
$$
This would contradict our choice of $\eta$ as he minimum diameter of the 
elements of~$\wlev{r}$.

Take $a\in A\setminus\bigcup\bigl\{B(\calB):\calB\in\mathfrak{H}\bigr\}$
and let $B$ be the closed $\epsilon$-ball about~$a$. 
We set $\varphi(\calA)=B\cap\bigcup\calA$.

If $E\in\calA$ then $\varrho(a,E)\le\varrho_H(A,E)<\epsilon$ and so $B$ 
meets~$E$; clearly $\diam B\le2\epsilon$, so $\varphi(\calA)$ has the first 
two required properties.
Finally, if $\calB\in\mathfrak{H}$ then 
$\varphi(\calA)\cap\varphi(\calB)\subseteq B\cap \varphi(\calB)=\emptyset$.

The collection
$$
\calN = \{\varphi(\calA) : \calA\in\mathfrak{F}\} 
$$
is finite, disjoint and has mesh less than $\half$.
To reach our final contradiction consider any continuum~$C$ in 
$X\setminus\bigcup\calN$, take $E\in\wlev{r}$ that intersects~$C$
and fix~$\calA\in\mathfrak{F}$ with $E\in\calA$. 

Now, $E$ meets $\varphi(\calA)$ and $C$ does not so $E\nsubseteq C$;
but then $C\subseteq E$, because $X$ is hereditarily indecomposable.
This, however, means that 
$$
    \diam C \le \diam E \le \half.
$$
This contradicts Lemma~\ref{lemma.small.mesh}.
\end{proof}

\begin{remark}
Although our argument took place in the hyperspace $\CsX$ it could have been
presented as a decomposition result as well.
We have already seen that $\wlev{r}$~is a decomposition of~$X$; because
$\wlev{r}$~is a closed subset of~$\CsX$ one can quite readily show that the 
decomposition map is actually closed and open.
The latter condition implies that the Hausdorff metric defines a compatible
metric on the decomposition space.
We find that $X$~admits an open continuous map onto an infinite-dimensional
continuum.
\end{remark}

\subsection{From three to infinity}

We now make good on our promise by showing that Brouwer's fixed-point
theorem for~$\I^3$ implies the full version.

As remarked before, Brouwer's Fixed-Point theorem for~$\I^n$ implies
that the pairs of faces of the cube~$\I^n$ form an essential family.
Thus, from the version for~$\I^3$ we find that the pairs $(A_0,B_0)$,
$(A_1,B_1)$ and $(A_2,B_2)$ form an essential family.
By Theorem~\ref{thm.hi.partition} let $L$ be any hereditarily indecomposable 
partition between~$A_2$ and~$B_2$.
By Lemma~\ref{lemma.partial.essential} the traces of the pairs~$(A_0,B_0)$,
and~$(A_1,B_1)$ form an essential family on~$L$, whence $\dim L\ge 2$.
There is a component~$X$ of~$S$ on which these traces also form an essential
family.
We find that Brouwer's theorem for~$\I^3$ implies the existence of an at least
two-dimensional hereditarily indecomposable continuum~$X$.
We shall now prove, from this fact, that for every~$m$ the cube~$\I^m$ has
the fixed-point property.

Working toward a contradiction we take the first~$m$ such that $\I^m$
has a fixed-point free map~$f$.
Using~$f$ one can make fixed-point free maps on every~$\I^k$ with~$k\ge m$
and, as is well-known, for every~$k\ge m$ a retraction of~$\I^k$ onto its 
boundary.

We use the Whitney level~$\wlev{r}$ and the number~$\eta$ from 
Proposition~\ref{prop.many.neighbours}.
To begin we set $n=3^{2m+1}-1$ and $\epsilon=\eta/(4n)$.
The compact space~$\wlev{r}$ has many finite open covers of mesh less 
than~$\epsilon$, each of which has a nerve, a polyhedron, associated with it,
see~\cite[1.10]{Engelking1995}.
The canonical map onto this nerve is an~$\epsilon$-map, i.e.,
each fiber has diameter less than~$\epsilon$.
We choose a polyhedron~$P$ of minimal dimension, say~$k$, such that
there is an $\epsilon$-map $f:\wlev{r}\to P$.

We use Proposition~\ref{prop.many.neighbours} to show that $k\ge m$.
Indeed, assume $k<m$ and apply Theorem~1.10.4 from~\cite{Engelking1995}
to see that $P$ may be realized inside~$\I^{2k+1}$ 
(here is where the Linear Algebra is needed).
For every~$l$ we can create a closed cover~$\calF_l$ of $\I^{2k+1}$ and hence 
of~$P$ by cutting along the hyperplanes $x_i=j/l$, where $i<2k+1$ and
$j=0$, \dots,~$l$.
Observe that every element of~$\calF_l$ meets at most $3^{2k+1}-1$ other 
elements of~$\calF_l$.
If $l$ is taken large enough then the preimage under~$f$ of~$\calF_l$
is a finite closed cover of~$\wlev{r}$ of mesh less than~$\epsilon$
such that every element meets fewer than~$n$ other elements of the cover.
This contradicts Proposition~\ref{prop.many.neighbours}.

We find that $k\ge m$.
To reach our final contradiction we consider the successive barycentric
subdivisions of~$P$.
In each of these subdivisions we find retractions of the $k$-simplices
onto their boundaries and combine these into a map~$r:P\to Q$,
where $Q$~is the union of the at most $(k-1)$-dimensional simplices
in the subdivision.
For a fine enough subdivision the composition $r\circ f$ is an
$\epsilon$-map from $\wlev{r}$ onto a $(k-1)$-dimensional polyhedron~$Q$.
This contradicts the minimality of~$k$.

\begin{remark}\label{rem.wlev(r).inf.dm}
The arguments given above imply in particular that $\wlev{r}$ cannot
be embedded into~$\I^n$ for any~$n$.
The Embedding Theorem (\cite[1.11.4]{Engelking1995}) now implies that
$\wlev{r}$~is infinite-dimensional.

This provides another route to Brouwer's Fixed-point theorem.
The first step is to observe that $\wlev{r}$ has arbitrarily large
finite essential families of pairs of closed sets.
The third step is to derive the fixed-point theorem for~$\I^n$ from the fact
that the faces of~$\I^n$ form an essential family, 
see \cite[1.8.B]{Engelking1995}.
The intermediate step is provided in the following proposition, which is
related to a theorem of Holszty\'nski from~\cite{Holsztynski1964}.
\end{remark}

\begin{proposition}
If some normal space~$X$ has an essential family consisting of $n$~pairs
then the pairs of opposite faces of\/~$\I^n$ also form an essential family.
\end{proposition}

\begin{proof}
Let $\bigl\{(C_i,D_i):i<n\bigr\}$ be an essential family in the normal 
space~$X$.
Apply Urysohn's lemma to get continuous functions $f_i:X\to\I$ such that
$f_i[C_i]=\{0\}$ and $f_i[D_i]=\{1\}$ for all~$i$ and take the diagonal 
map $f=\mapdiag_{i<n}f_i$.
If $L_i$ is a partition between the faces $A_i$ and $B_i$ of~$\I^n$
for each~$i$, then $f\preim[L_i]$ is a partition between $C_i$ and $D_i$
and so $\bigcap_{i<n}f\preim[L_i]\neq\emptyset$; but then 
$\bigcap_{i<n}L_i\neq\emptyset$ as well.
\end{proof}


\providecommand{\bysame}{\leavevmode\hbox to3em{\hrulefill}\thinspace}
\providecommand{\MR}{\relax\ifhmode\unskip\space\fi MR }
\providecommand{\MRhref}[2]{%
  \href{http://www.ams.org/mathscinet-getitem?mr=#1}{#2}
}
\providecommand{\href}[2]{#2}

\end{document}